\documentclass [12pt]  {article} 
\usepackage{amssymb}
\def \CC{{\mathbb{C}}}

\def \FFF{{\cal{F}}}

\begin{document}

\begin{center}
{\Large {\bf Entire functions sharing simple $a$-points with their 
first derivative}}\\
\bigskip
{\sc Andreas Schweizer\footnote{This paper was written while the 
author was working at the Institute of Mathematics of Academia Sinica
in Taipei, supported by grant 99-2115-M-001-011-MY2 from the 
National Science Council (NSC) of Taiwan.}}\\
\bigskip
{\small {\rm Department of Mathematics,\\
Korea Advanced Institute of Science and Technology (KAIST),\\ 
Daejeon 305-701\\
South Korea\\
e-mail: schweizer@kaist.ac.kr}}
\end{center}
\begin{abstract}
\noindent
We show that if a nonconstant complex entire function $f$ and 
its derivative $f'$ share their simple zeroes and their simple 
$a$-points for some nonzero constant $a$, then $f\equiv f'$. 
We also discuss how far these conditions can be relaxed or 
generalized. Finally, we determine all entire functions $f$ 
such that for $3$ distinct complex numbers $a_1,a_2,a_3$ 
every simple $a_j$-point of $f$ is an $a_j$-point of $f'$.
\end{abstract}

\subsection*{1. Introduction and results}

Throughout $f(z)$ or $f$ denotes an entire function, i.e., a function that 
is holomorphic in the whole complex plane, and $f'(z)$ or $f'$ denotes its 
derivative. We write $f\equiv g$ to say that the two functions are identical.
\par
Everybody knows that $f\equiv f'$ if and only if $f(z)=Ce^z$ with some 
complex constant $C$. For an apparently much weaker condition that has 
the same implication we recall the following.
\par
Two meromorphic functions $f$ and $g$ are said to share the value 
$a\in\CC$ IM (ignoring multiplicity), or just to share the value 
$a$, if $f$ takes the value $a$ at exactly the same points as $g$.
If moreover at any given such point the functions $f$ and $g$ take 
the value $a$ with the same multiplicity, then $f$ and $g$ are said 
to share the value $a$ CM (counting multiplicity).
\par
The general philosophy is that two meromorphic functions that share
``too many'' values must be equal. In the special situation where
$g$ is the derivative of $f$ one usually needs even less values.
\\ \\
{\bf Theorem A.} [RuYa] \it
Let $f$ be a nonconstant entire function and let $a$ and $b$ be complex 
numbers with $a\neq b$. If $f$ and $f'$ share the values $a$ and $b$ CM, 
then $f\equiv f'$.
\rm
\\ \\
Actually, in this case one doesn't even need the multiplicities, as was 
proved two years later.
\\ \\
{\bf Theorem B.} [MuSt, Satz 1] \it
Let $f$ be a nonconstant entire function and let $a$ and $b$ be complex 
numbers with $a\neq b$. If $f$ and $f'$ share the values $a$ and $b$ IM, 
then $f\equiv f'$.
\rm
\\ \\
See also [YaYi, Theorem 8.3] for a proof in English. Theorem A and 
later Theorem B have been generalized in [LiYi] resp. [L\"uXuYi] 
by relaxing the requirements on the sharing. We content ourselves 
with one representative example, which we will need later.
\\ \\
{\bf Theorem C.} [L\"uXuYi, Corollary 1.1] \it
Let $b$ be a nonzero number and let $f(z)$ be a nonconstant entire function.
If $f(z)=0\Rightarrow f'(z)=0$ and $f(z)=b\Rightarrow f'(z)=b$, then one of
the following cases must occur:
\begin{itemize}
\item[(a)] $f\equiv f'$,
\item[(b)] $f=b\{\frac{1}{4}A^2e^{z/2}+Ae^{z/4}+1\}$, where $A$ is a nonzero 
constant.
\end{itemize}
\rm
\bigskip
\noindent
There are yet other ways to generalize the sharing of two values between
$f$ and $f'$. For example the papers [LiYa] and [L\"uXu] also treat the 
case of $f$ and $f'$ sharing a two-element-set $\{a,b\}$, i.e., the two 
shared values might get mixed up.
\par
But we want to investigate a generalization that, to the best of our 
knowledge, has not been applied yet to an entire function and its 
derivative.
\\ \\
{\bf Definition.}
Let $f$ and $g$ be two meromorphic functions and $a\in\CC$. We say that
$f$ and $g$ {\bf share their simple $a$-points} if the points where $f$ 
takes the value $a$ with multiplicity one are exactly those where $g$ 
takes the value $a$ with multiplicity one.
\par
In the literature sometimes the notation 
$$\overline{E}_{1)}(a,f)=\overline{E}_{1)}(a,g)$$
is used to describe this kind of sharing.
\\ \\
Obviously, this property generalizes sharing $a$ CM, and even sharing $a$ 
with weight one in the sense of [La], but in general it is neither stronger 
nor weaker than sharing $a$ IM. Note that sharing simple $a$-points does 
{\it not} imply sharing the value $a$ since we make no requirements at all
concerning points where the value $a$ is taken with higher multiplicity.
\par
Very roughly, the philosophy is that the bulk of $a$-points should be simple,
at least if one considers enough values $a$, and hence the loss of information 
when giving up control over the multiple values can be compensated.
\par
For example, in contrast to Nevanlinna's famous theorem that two nonconstant 
meromorphic functions that share $5$ values must be equal, one obtains that 
two nonconstant meromorphic functions that share simple $a_j$-points for 
$7$ values $a_j$ must be equal [GoBo] or [YaYi, Section 3.3.1]. 
\par
In this paper we try to find similar results corresponding to 
Theorems A and B.
\\ \\
{\bf Theorem 1.} \it
Let $f(z)$ be a nonconstant entire function and $0\neq a\in\CC$. If $f$ and
its derivative $f'$ share their simple $a$-points and their simple zeroes, 
then $f\equiv f'$. 
\\ \\
\rm 
Actually, we even prove a slightly stronger statement in the spirit of
[LiYi] and [L\"uXuYi].
We will be working with the following conditions (in order of decreasing 
strength):
\begin{itemize}
\item $f$ and $f'$ share their simple $a$-points, i.e.
$$(f=a\ \hbox{\rm and}\ f'\neq 0)\ \Leftrightarrow 
\ (f'=a\ \hbox{\rm and}\ f''\neq 0);$$
\item Every simple $a$-point of $f$ is a simple $a$-point of $f'$, i.e.
$$(f=a\ \hbox{\rm and}\ f'\neq 0)\ \Rightarrow
\ (f'=a\ \hbox{\rm and}\ f''\neq 0);$$
\item Every simple $a$-point of $f$ is a (not necessarily simple) $a$-point 
of $f'$, i.e.
$$(f=a\ \hbox{\rm and}\ f'\neq 0)\ \Rightarrow\ f'=a.$$
This condition is of course equivalent to
$$f=a\Rightarrow f'\in\{a,0\}.$$
\end{itemize}
\bigskip
\noindent
{\bf Theorem 1'.} \it
Let $f(z)$ be a nonconstant entire function and $0\neq a\in\CC$. 
If $f$ and $f'$ share their simple zeroes and if every simple 
$a$-point of $f$ is a (not necessarily simple) $a$-point of $f'$,
then $f\equiv f'$. 
\\ \\
\rm
The proof is given in Section 2.
\\ \\
{\bf Example 1.}
From [LiYi, Theorem 2] we take the function 
$$f(z)=Ce^{\frac{b}{b-a}z}+a$$ 
with nonzero constants $C$, $a$, $b(\neq a)$.
It shares the value $b$ CM with $f'$ and omits the value $a$.
This shows that in Theorem 1' we cannot simply replace sharing simple
zeroes by sharing simple $b$-points. 
This is perhaps not overly surprising. As the easy Lemma 1 in the next 
section shows, sharing simple zeroes with $f'$ has much stronger
implications on $f$ than sharing simple $b$-points for some $b\neq 0$.
\\ \\
{\bf Example 2.}
Let 
$$f(z)=\frac{a}{2}(\sin(2z)+1);$$
then $f'(z)=a\cos(2z)$. All $a$-points of $f$ and of $f'$ and all zeroes
of $f$ have multiplicity $2$. Thus the condition that $f$ and $f'$ share
their simple $a$-points and that every simple zero of $f$ is a simple zero 
of $f'$ does not imply $f\equiv f'$.
\\ \\
The more interesting question is whether in Theorem 1 we can replace 
sharing the simple zeroes by sharing the simple $b$-points for some 
nonzero $b$ different from $a$.
In general the answer again is negative.
\\ \\
{\bf Example 3.}
Let $0\neq a\in\CC$. The entire function $f(z)= a\sin z$ and its derivative
 $f'(z)=a\cos z$ share their simple $a$-points and their simple $-a$-points,
for the trivial reason that all their $a$-points and $-a$-points have 
multiplicity $2$.
\\ \\
However, somehow this counterexample seems to hinge on the fact that the
second value is the negative of the first. It is still conceivable that
if $f$ and $f'$ share their simple $a$-points and their simple $b$-points
the sufficient condition that forces $f\equiv f'$ is simply $a+b\neq 0$,
not $ab=0$. 
\par
For example, in the somewhat similar context of $f$ and $f'$ sharing 
a two-element-set $\{a,b\}$ CM, the only case for which non-obvious 
functions $f$ exist is $a+b=0$ (see [LiYa, Theorem 3] and [L\"uXu]).
\par
So we ask the following
\\ \\
{\bf Question.}
Let $a$, $b$ be two distinct nonzero complex numbers with $a+b\neq 0$. 
If a nonconstant entire function $f$ and its derivative $f'$ share 
their simple $a$-points and their simple $b$-points, does this imply 
$f\equiv f'$?
\\ \\
At the moment we don't know the answer and we do not even have a clear 
feeling whether it will be positive or negative.
As a small consolation we prove another result, which in case of a positive
answer to this question would follow as an immediate corollary.
\\ \\
{\bf Theorem 2.} \it
Let $a_1$, $a_2$, $a_3$ be three distinct complex numbers and let 
$f$ be a nonconstant entire function. If $f$ and $f'$ share their 
simple $a_j$-points for $j=1,2,3$, then $f\equiv f'$.
\\ \\
\rm
Again, we prove a stronger result.
\\ \\
{\bf Theorem 2'.} \it
Let $a_1$, $a_2$, $a_3$ be three distinct complex numbers. Nonconstant 
entire functions $f$ with $f\not\equiv f'$ and 
$$f=a_j\Rightarrow f'\in\{a_j,0\}$$
for $j=1,2,3$ exist if and only if $a_j=\zeta^ja_3$ with
$\zeta$ being a third root of unity, that is, if $(X-a_1)(X-a_2)(X-a_3)$
is of the form $X^3-\delta$.
\par
Moreover, functions with this property necessarily are of the form
$$f(z)=\frac{4\delta}{27\beta^2}e^{\frac{2}{3}z}+\beta e^{-\frac{1}{3}z}$$
with a nonzero constant $\beta$.
\par
Conversely, every function of this form has the stronger property that
every simple $a_j$-point of $f$ is a simple $a_j$-point of $f'$ for
$j=1,2,3$.
\\ \\
\rm
The three Examples above and the second case of Theorem C show that the 
condition in Theorem 2' cannot be reduced to two values $a_1$, $a_2$.
\\

\subsection*{2. Proofs}

The following observation is almost trivial.
\\ \\
{\bf Lemma 1.} \it
Suppose that $f$ and $f'$ share their simple zeroes. Then
\begin{itemize}
\item[(a)] $f'$ has no simple zeroes.
\item[(b)] Every multiple point of $f$ has multiplicity at least $3$.
\item[(c)] Every zero of $f$ has multiplicity at least $3$.
\end{itemize}
\rm
\bigskip
\noindent
The proofs of the theorems follow an overall strategy that we have seen 
in several articles from the last ten years on entire (or meromorphic) 
functions $f$ with certain value sharing properties. This strategy gains 
its strength from the combination of different methods. To emphasize 
this we have divided it into four steps.
\par
In Step 1 one constructs from $f$ a family of analytic functions
by shifting the argument and then uses the properties of $f$ to 
show that this family is normal.
\par
In Step 2 one obtains from the normality of that family that $f$ has 
order at most $1$. This is almost automatic. Nevertheless, we want 
to consider this as a separate step for the following reason: Even
if the family in Step 1 is not normal, there might be other ways
to show that $f$ has order at most $1$.
\par
In Step 3 one constructs an auxiliary function $h$ from $f$ and its 
derivative(s) that somehow encodes the value sharing property of $f$.
Then one uses Nevanlinna Theory arguments to show that $h$ is constant.
This task is greatly facilitated, and sometimes only possible, thanks
to the knowledge that $f$ has order at most $1$.
\par
In Step 4 one uses the information encoded in $h\equiv const$ to derive
the desired properties of $f$. In the proof of Theorem 2' we will to 
that end emphasize geometric considerations concerning the algebraic 
curve described by $h\equiv const$.
\par
To start with, i.e. for Step 1, we need the following minor strengthening
of the famous Zalcman Lemma.
\\ \\
{\bf Lemma 2.} \it 
Let $\FFF$ be a family of holomorphic functions on the unit disk.
If $\FFF$ is not normal, then there exist
\begin{itemize}
\item[(i)] a number $0<r<1$,
\item[(ii)] points $z_n$, $|z_n|<r$,
\item[(iii)] functions $f_n\in\FFF$,
\item[(iv)] positive numbers $\rho_n\to 0$,
\end{itemize}
such that
$$f_n(z_n +\rho_n\xi)=:g_n(\xi)\to g(\xi)$$
uniformly on compact subsets of $\CC$, where $g$ is a nonconstant 
entire function.
\par
Moreover, given a complex number $a$, if there exists a bound 
$M$ and a positive integer $m$ such that for every function 
$f$ in $\FFF$ and every $z_0\in\CC$ with $f(z_0)=a$ we have
$|f^{(k)}(z_0)|\le M$ for $k=1,2,\ldots m$,
then every $a$-point of $g$ has multiplicity at least $m+1$.
\\ \\
{\bf Proof.} \rm\ 
This is essentially the original version of Zalcman's Lemma
[Za]. The only thing we have to prove is the last assertion.
For fixed $n$ we differentiate $g_n(\xi)$ with respect to 
$\xi$ and get
$$g_n^{(k)}(\xi)=\rho_n^kf_n^{(k)}(z_n+\rho_n\xi).$$
Now suppose that $g(\xi_0)=a$. Since $g$ is nonconstant, by
Hurwitz's theorem there exist $\xi_n$, $\xi_n\to\xi_0$, such
that for sufficiently large $n$ we have
$a=g(\xi_0)=g_n(\xi_n)$. Hence our assumptions imply
$|g_n^{(k)}(\xi_n)|\leq\rho_n^k M$ for $k=1,2,\ldots,m$ and 
$n$ sufficiently large. Since $g_n^{(k)}(\xi)$ converges 
locally uniformly to $g^{(k)}(\xi)$, we obtain
$$g^{(k)}(\xi_0)=\lim_{n\to\infty}g_n^{(k)}(\xi_n)=0$$
for $k=1,2,\ldots,m.$
\hfill$\Box$
\\ \\
{\bf Remark.} \rm
For $m=1$ the same argument is already in [Ch] (and in other
papers). 
\par
Actually, [Ch] claims that under certain stronger conditions 
the function $g$ omits the value $a$. However, we feel uneasy 
about Theorem 1 and Theorem 3 in [Ch] as it seems to us that 
by practically the same argument one would then be able to 
prove that a holomorphic family $\FFF$ with 
$f(z)\in\{1,-1\}\Rightarrow f'(z)\in\{1,-1\}$ for all $f\in\FFF$
would be normal, in contradiction to Example 1 in the same paper.
\par
The problem seems to be that on lines 6 and 7 of page 1476 of 
[Ch] the value $a_l$ in $g_n'(\xi_n^{(j)})=\rho_n a_l$ depends 
on $j$, and therefore on line 13 one cannot conclude that 
$g_n'(\xi)-\rho_n a_l$ has $k$ zeroes.
\\ \\
{\bf Proof of Theorem 1'.} \rm\ 
\\ \\
{\bf Step 1:} We consider the family of holomorphic functions
$\FFF=\{f_\omega(z)\ :\ \omega\in\CC\}$ with 
$f_\omega(z)=f(z+\omega)$. Note that due to its special form 
$\FFF$ is normal on $\CC$ if and only if it is normal
on the unit disk. Obviously, this family satisfies the
conditions of Lemma 2 for $0$ with $m=2$ and for $a$ with $m=1$.
We conclude that $\FFF$ must be normal. Otherwise we could 
construct a nonconstant entire function $g$ such that all 
$a$-points of $g$ have multiplicity at least $2$ and all
zeroes have multiplicity at least $3$. So for the function 
$\Theta$ (the sum of the deficiency and the ramification 
defect) we would have $\Theta(a,g)\geq\frac{1}{2}$ and 
$\Theta(0,g)\geq\frac{2}{3}$, in contradiction to the defect 
relation $\sum_{b\in\CC}\Theta(b,g)\leq 1$ for entire functions
([ChYe, Corollary 5.2.4] or [YaYi, Section 1.2.4]).
\\ \\
{\bf Step 2:} From Step 1 we readily obtain that $f$ has order at most 
$1$. This is a general principle; $f$ is a Yosida function (i.e. its 
spherical derivative is uniformly bounded on $\CC$) if and only if the 
family $\{f(z+\omega)\ :\ \omega\in\CC\}$ is normal on $\CC$ [Mi, p.198], 
and a holomorphic Yosida function has order at most $1$ [Mi, p.211].
\\ \\
{\bf Step 3:} Consider the auxiliary function
$$h=\frac{(f')^2(f-f')}{f^2(f-a)}.$$
It is easy to see that the potential poles arising from zeroes of 
$f-a$ are cancelled either by $f'=a$ or by the zeroes of $(f')^2$.
As for the zeroes of $f$, note that by Lemma 1(c) then $f-f'$ has 
at least a double zero. So $h$ is an entire function.
\par
Using the standard functions from Nevanlinna theory and their basic 
properties (see e.g. [ChYe] or [YaYi]), from
$$h=\frac{f'}{f}\cdot\frac{f'}{f-a}-\frac{f'}{f}\cdot\frac{f'}{f}\cdot
\frac{f'}{f-a}$$
we obtain
\begin{eqnarray*}
T(r,h) & = & m(r,h)\\
       &\leq & m(r,\frac{f'}{f})+m(r,\frac{f'}{f-a})
+m(r,\frac{f'}{f})+m(r,\frac{f'}{f})+m(r,\frac{f'}{f-a})+O(1).
\end{eqnarray*}
From [NgOs] or [HeKoR\"a, Theorem 4.1] we see that if $f$ is an entire 
function of order at most $1$, then $m(r,\frac{f'}{f})=o(\log r)$
(compare [ChYe, Section 3.5]). Hence the above estimate gives 
$T(r,h)=o(\log r)$, which means that $h$ is constant.
\\ \\
{\bf Step 4:} If $h\equiv 0$, then $(f')^2(f-f')\equiv 0$, and hence
$f\equiv f'$ since $f$ is nonconstant.
\par
Now consider the case $h\equiv\gamma$ for some nonzero constant $\gamma$.
Then every $a$-point of $f$ must be simple; otherwise by Lemma 1(b)
it would have multiplicity at least $3$ and then the term $(f')^2$
would cause a zero of $h$. So we have $f=a\Rightarrow f'=a$ and by
Lemma 1(c) also $f=0\Rightarrow f'=0$, that is, we are in the situation
of Theorem C. But the second possibility 
$f=a(\frac{A}{2}e^{\frac{z}{4}}+1)^2$
is ruled out, for example because all zeroes of that function have
multiplicity $2$, in contradiction to Lemma 1(c).
\hfill $\Box$
\\ \\
{\bf Proof of Theorem 2'.} \rm\ \\
We prefer to write $a,b,c$ for $a_1,a_2,a_3$.
\par
The holomorphic family $\{f_\omega(z)\ :\ \omega\in\CC\}$ with 
$f_\omega(z)=f(z+\omega)$ is normal. If not, as in Step 1 of the proof 
of Theorem 1' we could construct a nonconstant entire function $g$ 
with $3$ totally ramified values (namely $a$, $b$ and $c$); 
this would contradict the defect relations [ChYe, Theorem 5.4.1].
But actually our claim is just a special case of [FaZa, Lemma 4].
\par
Next, exactly the same argument as in Step 2 shows that $f$ has order
at most $1$, and as in Step 3 we see that
$$h=\frac{(f')^2(f-f')}{(f-a)(f-b)(f-c)}$$
$$=\frac{f'}{f-c}\left(\frac{b}{b-a}\cdot\frac{f'}{f-b}-\frac{a}{b-a}\cdot
\frac{f'}{f-a}-\frac{f'}{f-a}\cdot\frac{f'}{f-b}\right)$$
must be constant. 
If $h\equiv 0$, again we get $f\equiv f'$ since $f'\not\equiv 0$. 
\par
Now we discuss the case $h\equiv\gamma$ for a nonzero constant $\gamma$.
This can probably be done by some case distinctions as in [L\"uXuYi], 
or rather, more complicated ones. But we prefer a more geometric argument 
that would also work in many more complicated situations.
\par
Consider the holomorphic map 
$$z\mapsto(f(z),f'(z))=(X,Y)$$ 
from $\CC$ to the affine curve
$$A:\ \ Y^3-XY^2+\gamma(X-a)(X-b)(X-c)=0.$$
This polynomial is irreducible in $\CC[X,Y]$. If not, it would have 
a factor $Y-uX+v$ with $u\neq 0$; but then by the equation above
$X-\frac{v}{u}$ would be a multiple factor of $(X-a)(X-b)(X-c)$.
\par
Now we write $(X-a)(X-b)(X-c)$ as $X^3+c_2X^2+c_1X+c_0$ and examine
the corresponding projective curve
$$R:\ \ Y^3-XY^2+\gamma(X^3+c_2X^2Z+c_1XZ^2+c_0Z^3)=0.$$
This is an irreducible cubic curve. So either it is smooth and has 
genus $1$, or it has exactly one singular point and genus $0$. In 
the latter case the smooth model of $R$ is a Riemann sphere.
We suppress discussing the somewhat complicated conditions on 
$\gamma$ and $c_2, c_1,c_0$ that distinguish the two cases, as it
would not really help us in finishing the proof.
\par
If the genus is $1$, the affine curve $A$ is obtained by removing
at least one point from a smooth, projective curve (equivalently,
from a compact Riemann surface) of positive genus. Hence $A$ is 
hyperbolic [Fo, Theorem 27.12], i.e., its universal covering is the 
unit disk. Hence (essentially by Liouville's Theorem) every
holomorphic map from $\CC$ to $A$ must be constant. This would 
mean that $f$ is constant. 
\par
Alternatively, as Andreas Sauer has pointed out to me, hyperbolicity
of algebraic curves can also be obtained from value distribution 
theory of meromorphic functions. See [Ne, Chapter X, \S 3] and the 
references given there. As mentioned there, it is already a classical 
theorem by Picard [Pi] that if an algebraic curve $F(X,Y)=0$ is 
uniformized by two nonconstant meromorphic functions $X(z)$ and $Y(z)$
then the curve necessarily has genus $0$ or $1$. But if the genus is 
$1$, the functions are elliptic and hence not entire.
\par
Either way, we can assume from now on that the genus of $R$ is $0$. 
Then the function field $\CC(f,f')$ is a rational function field 
$\CC(t)$. It is a classical result that $f'$ has the same order 
(in the sense of Nevanlinna theory) as $f$ and that adding, 
multiplying and dividing functions does not increase the order.
As some textbooks do not mention this, we give the reference
[YaYi, Theorem 1.21 and Section 1.3.4]. 
Since $t$ is a rational expression in $f$ and $f'$, we thus obtain 
that $t(z)$ is a meromorphic function of order at most $1$.
\par
Plugging $Z=0$ into the homogeneous equation for $R$, we see that 
there are at least $2$ points outside the affine part. After a
M\"obius transformation we can assume that $t$ has a zero and 
a pole at these two points. Then $t(z)$ is an entire function
of order at most $1$ without zeroes. By Hadamard's factorization 
theorem [YaYi, Theorem 2.5] we have
$$t(z)=\frac{e^{\alpha z}}{\beta}$$
with nonzero constants $\alpha,\beta$.
\par
Of course, $f$ is a rational function of degree $3$ in $t$. If 
$t(z)$ omits other values than $\infty$ and $0$, it must be constant 
by Picard's Theorem, so $f$ would be constant. 
Thus the poles of $f$ are exactly at $t=0$ and $t=\infty$. Replacing 
$t$ by $\frac{1}{t}$ if necessary, we can assume that the double pole 
is at $t=\infty$. Then $f=\frac{b_2t^3+b_1t^2+b_0t+b_{-1}}{t}$ with 
$b_2b_{-1}\neq 0$. Choosing $\beta$ suitably, we can assume $b_{-1}=1$,
that is,
$$f=b_2t^2+b_1t+b_0+\frac{1}{t},$$
and hence
$$f'=2\alpha b_2t^2+\alpha b_1 t-\frac{\alpha}{t}.$$
We plug this into the equation for $A$ and compare coefficients for
the powers of $t$. From the coefficients of $t^6$ and of $t^{-3}$ we 
obtain 
$8\alpha^3 b_2^3 -4\alpha^2 b_2^3 +\gamma b_2^3 =0$ and
$-\alpha^3 -\alpha^2 +\gamma=0$, so together
$$\alpha=\frac{1}{3}\ \ \hbox{\rm and}\ \ \gamma=\frac{4}{27}.$$
Using this, the coefficient of $t^4$ forces $c_2=0$, and then from 
the $t^{-2}$-coefficient $b_0 =0$ follows. From the $t^3$-coefficient
we get $b_1 =0$, and with that the $t^{-1}$-coefficient implies 
$c_1 =0$. Finally, the coefficient of $t^0$ tells us that 
$b_2 =\frac{-4}{27}c_0 =\frac{4}{27}\delta$. This shows
$$f(z)=\frac{4\delta}{27\beta^2}e^{\frac{2}{3}z}+\beta e^{-\frac{1}{3}z}$$
and proves the main part of the theorem.
\\ \\
Checking the coefficients for the remaining powers of $t$ confirms
that such $f$ do indeed satisfy the differential equation
$$(f')^3-f(f')^2+\frac{4}{27}(f^3-\delta)\equiv 0.$$
So, conversely, assume $f$ is of the above form. Then the differential 
equation shows that if $f=\zeta^j c$ and $f'\neq 0$ then $f'=\zeta^j c$.
Moreover, since $f$ also satisfies the differential equation
$f''\equiv\frac{1}{3}f'+\frac{2}{9}f$,
such $\zeta^jc$-points of $f'$ are simple.
\hfill $\Box$
\\ \\
{\bf Proof of Theorem 2.} \rm\ \\
If $f\not\equiv f'$, then by Theorem 2' we must have $a_j=\zeta^j c$
and $f=\frac{4c^3}{27}t^2 +\frac{1}{t}$ with 
$t=\frac{1}{\beta}e^{\frac{1}{3}z}$. Correspondingly, 
$f'=\frac{8c^3}{81}t^2 -\frac{1}{3t}$ and
$f''=\frac{16c^3}{243}t^2 +\frac{1}{9t}$. In particular, 
$f'=c$ if and only if 
$t\in\{\frac{-3}{c}, \frac{3(2\pm\sqrt{6})}{4c}\}$.
But for $t=\frac{3(2+\sqrt{6})}{4c}$ we get $f''\neq 0$ and
$f=(\sqrt{6}-\frac{1}{2})c\neq c$. So not every simple $c$-point 
of $f'$ is a $c$-point of $f$.
\hfill$\Box$
\\

\subsection*{\hspace*{10.5em} References}
\begin{itemize}

\item[{[Ch]}] J.-F. Chen: \rm Shared sets and normal families 
of meromorphic functions, \it Rocky Mountain J. Math. 
\bf vol. 40, no. 5 \rm (2010), 1473-1479

\item[{[ChYe]}] W. Cherry and Z. Ye: \it Nevanlinna's Theory of
Value Distribution, \rm Springer, Berlin Heidelberg New York, 2001

\item[{[FaZa]}] M. Fang and L. Zalcman: \rm Normal families and 
uniqueness theorems for entire functions,  \it J. Math. Anal. Appl. 
\bf 280  \rm (2003), 273-283 

\item[{[Fo]}] O. Forster: \it Lectures on Riemann Surfaces,
\rm Springer Graduate Texts in Mathematics 81, 
Berlin Heidelberg 
New York, 1981

\item[{[GoBo]}] H.S. Gopalakrishna and S.S. Bhoosnurmath: \rm 
Uniqueness theorems for meromorphic functions, \it Math. Scand.
\bf 39 \rm (1976), 125-130

\item[{[HeKoR\"a]}] J. Heittokangas, R. Korhonen and J. R\"atty\"a:
\rm Generalized logarithmic derivative estimates of Gol'dberg-Grinshtein
type, \it Bull. London Math. Soc. \bf 36 \rm (2004), 105-114

\item[{[La]}] I. Lahiri: \rm Weighted sharing and uniqueness of
meromorphic functions, \it Nagoya Math. J. \bf 161 \rm (2001), 193-206

\item[{[LiYa]}] P. Li and C.-C. Yang: \rm Value sharing of an entire 
function and its derivatives, \it J. Math. Soc. Japan \bf 51 no. 4
\rm (1999), 781-799

\item[{[LiYi]}] J. Li and H. Yi: \rm Normal families and uniqueness
of entire functions and their derivatives, \it Arch. Math. (Basel)
\bf 87 \rm (2006), 52-59

\item[{[L\"uXu]}] F. L\"u and J. Xu: \rm Sharing set and normal 
families of entire functions and their derivatives, \it Houston J. Math.
\bf vol. 34, no. 4 \rm (2008), 1213-1223

\item[{[L\"uXuYi]}] F. L\"u, J. Xu and H. Yi: \rm Uniqueness 
theorems and normal families of entire functions and their derivatives, 
\it Ann. Polon. Math. \bf 95.1 \rm (2009), 67-75

\item[{[Mi]}] D. Minda: \rm Yosida functions, in: \it Lectures on
Complex Analysis, Xian 1987, (C.T. Chuang, ed.) \rm World Scientific,
Singapore, 1988, pp.197-213

\item[{[MuSt]}] E. Mues and N. Steinmetz: \rm Meromorphe Funktionen, 
die mit ihrer Ableitung Werte teilen, \it Manuscripta Math. \bf 29 \rm 
(1979), 195-206

\item[{[Ne]}] R. Nevanlinna: \it Eindeutige analytische Funktionen,
\rm (2nd edition), Springer Grundlehren vol. 46, 
Berlin G\"ottingen Heidelberg, 1953

\item[{[NgOs]}] V. Ngoan and I.V. Ostrovskii: \rm The logarithmic 
derivative of a meromorphic function, \it Akad. Nauk. Armyan. SSR Dokl. 
\bf 41 \rm (1965), 272-277 (in Russian)

\item[{[Pi]}] E. Picard: \rm D\'emonstration d'un th\'eor\`eme 
g\'en\'eral sur les fonctions uniformes li\'ees par une relation 
alg\'ebrique, \it Acta Math. \bf 11 \rm (1887), 1-12 

\item[{[RuYa]}] L. Rubel and C. C. Yang: \rm Values shared by an entire
function and its derivative, in: \it Complex Analysis, Kentucky 1976, 
(J.D. Buckholtz and T.J. Suffridge, eds.)
\rm Springer LNM 599, Berlin, 1977, pp.101-103

\item[{[YaYi]}] C.-C. Yang and H.-X. Yi: \it Uniqueness theory
of meromorphic functions, \rm Kluwer Academic Publishers Group,
Dordrecht, 2003

\item[{[Za]}] L. Zalcman: \rm A Heuristic Principle in Complex 
Function Theory, \it Amer. Math. Monthly \bf 82 \rm (1975), 813-817

\end{itemize}

\end{document}